
\ifx\shlhetal\undefinedcontrolsequence\let\shlhetal\relax\fi

\input amstex
\expandafter\ifx\csname mathdefs.tex\endcsname\relax
  \expandafter\gdef\csname mathdefs.tex\endcsname{}
\else \message{Hey!  Apparently you were trying to
  \string\input{mathdefs.tex} twice.   This does not make sense.} 
\errmessage{Please edit your file (probably \jobname.tex) and remove
any duplicate ``\string\input'' lines}\endinput\fi




\catcode`\X=12\catcode`\@=11

\def\n@wcount{\alloc@0\count\countdef\insc@unt}
\def\n@wwrite{\alloc@7\write\chardef\sixt@@n}
\def\n@wread{\alloc@6\read\chardef\sixt@@n}
\def\r@s@t{\relax}\def\v@idline{\par}\def\@mputate#1/{#1}
\def\l@c@l#1X{\firstpart.#1}\def\gl@b@l#1X{#1}\def\t@d@l#1X{{}}

\def\crossrefs#1{\ifx\all#1\let\tr@ce=\all\else\def\tr@ce{#1,}\fi
   \n@wwrite\cit@tionsout\openout\cit@tionsout=\jobname.cit 
   \write\cit@tionsout{\tr@ce}\expandafter\setfl@gs\tr@ce,}
\def\setfl@gs#1,{\def\@{#1}\ifx\@\empty\let\next=\relax
   \else\let\next=\setfl@gs\expandafter\xdef
   \csname#1tr@cetrue\endcsname{}\fi\next}
\def\m@ketag#1#2{\expandafter\n@wcount\csname#2tagno\endcsname
     \csname#2tagno\endcsname=0\let\tail=\all\xdef\all{\tail#2,}
   \ifx#1\l@c@l\let\tail=\r@s@t\xdef\r@s@t{\csname#2tagno\endcsname=0\tail}\fi
   \expandafter\gdef\csname#2cite\endcsname##1{\expandafter
     \ifx\csname#2tag##1\endcsname\relax?\else\csname#2tag##1\endcsname\fi
     \expandafter\ifx\csname#2tr@cetrue\endcsname\relax\else
     \write\cit@tionsout{#2tag ##1 cited on page \folio.}\fi}
   \expandafter\gdef\csname#2page\endcsname##1{\expandafter
     \ifx\csname#2page##1\endcsname\relax?\else\csname#2page##1\endcsname\fi
     \expandafter\ifx\csname#2tr@cetrue\endcsname\relax\else
     \write\cit@tionsout{#2tag ##1 cited on page \folio.}\fi}
   \expandafter\gdef\csname#2tag\endcsname##1{\expandafter
      \ifx\csname#2check##1\endcsname\relax
      \expandafter\xdef\csname#2check##1\endcsname{}%
      \else\immediate\write16{Warning: #2tag ##1 used more than once.}\fi
      \multit@g{#1}{#2}##1/X%
      \write\t@gsout{#2tag ##1 assigned number \csname#2tag##1\endcsname\space
      on page \number\count0.}%
   \csname#2tag##1\endcsname}}

\def\multit@g#1#2#3/#4X{\def\t@mp{#4}\ifx\t@mp\empty%
      \global\advance\csname#2tagno\endcsname by 1 
      \expandafter\xdef\csname#2tag#3\endcsname
      {#1\number\csname#2tagno\endcsnameX}%
   \else\expandafter\ifx\csname#2last#3\endcsname\relax
      \expandafter\n@wcount\csname#2last#3\endcsname
      \global\advance\csname#2tagno\endcsname by 1 
      \expandafter\xdef\csname#2tag#3\endcsname
      {#1\number\csname#2tagno\endcsnameX}
      \write\t@gsout{#2tag #3 assigned number \csname#2tag#3\endcsname\space
      on page \number\count0.}\fi
   \global\advance\csname#2last#3\endcsname by 1
   \def\t@mp{\expandafter\xdef\csname#2tag#3/}%
   \expandafter\t@mp\@mputate#4\endcsname
   {\csname#2tag#3\endcsname\lastpart{\csname#2last#3\endcsname}}\fi}
\def\t@gs#1{\def\all{}\m@ketag#1e\m@ketag#1s\m@ketag\t@d@l p
\let\realscite\scite
\let\realstag\stag
   \m@ketag\gl@b@l r \n@wread\t@gsin
   \openin\t@gsin=\jobname.tgs \re@der \closein\t@gsin
   \n@wwrite\t@gsout\openout\t@gsout=\jobname.tgs }
\outer\def\localtags{\t@gs\l@c@l}
\outer\def\globaltags{\t@gs\gl@b@l}
\outer\def\newlocaltag#1{\m@ketag\l@c@l{#1}}
\outer\def\newglobaltag#1{\m@ketag\gl@b@l{#1}}

\newif\ifpr@ 
\def\m@kecs #1tag #2 assigned number #3 on page #4.%
   {\expandafter\gdef\csname#1tag#2\endcsname{#3}
   \expandafter\gdef\csname#1page#2\endcsname{#4}
   \ifpr@\expandafter\xdef\csname#1check#2\endcsname{}\fi}
\def\re@der{\ifeof\t@gsin\let\next=\relax\else
   \read\t@gsin to\t@gline\ifx\t@gline\v@idline\else
   \expandafter\m@kecs \t@gline\fi\let \next=\re@der\fi\next}
\def\pretags#1{\pr@true\pret@gs#1,,}
\def\pret@gs#1,{\def\@{#1}\ifx\@\empty\let\n@xtfile=\relax
   \else\let\n@xtfile=\pret@gs \openin\t@gsin=#1.tgs \message{#1} \re@der 
   \closein\t@gsin\fi \n@xtfile}

\newcount\sectno\sectno=0\newcount\subsectno\subsectno=0
\newif\ifultr@local \def\ultralocal{\ultr@localtrue}
\def\firstpart{\number\sectno}
\def\lastpart#1{\ifcase#1 \or a\or b\or c\or d\or e\or f\or g\or h\or 
   i\or k\or l\or m\or n\or o\or p\or q\or r\or s\or t\or u\or v\or w\or 
   x\or y\or z \fi}

\def\resetall{\global\advance\sectno by 1\subsectno=0
   \gdef\firstpart{\number\sectno}\r@s@t}
\def\resetsub{\global\advance\subsectno by 1
   \gdef\firstpart{\number\sectno.\number\subsectno}\r@s@t}
\def\newsection#1\par{\resetall\vskip0pt plus.3\vsize\penalty-250
   \vskip0pt plus-.3\vsize\bigskip\bigskip
   \message{#1}\leftline{\bf#1}\nobreak\bigskip}
\def\subsection#1\par{\ifultr@local\resetsub\fi
   \vskip0pt plus.2\vsize\penalty-250\vskip0pt plus-.2\vsize
   \bigskip\smallskip\message{#1}\leftline{\bf#1}\nobreak\medskip}


\newdimen\marginshift

\newdimen\margindelta
\newdimen\marginmax
\newdimen\marginmin

\def\margininit{       
\marginmax=3 true cm                  
				      
\margindelta=0.1 true cm              
\marginmin=0.1true cm                 
\marginshift=\marginmin
}    

\def\t@gsjj#1,{\def\@{#1}\ifx\@\empty\let\next=\relax\else\let\next=\t@gsjj
   \def\@@{p}\ifx\@\@@\else
   \expandafter\gdef\csname#1cite\endcsname##1{\citejj{##1}}
   \expandafter\gdef\csname#1page\endcsname##1{?}
   \expandafter\gdef\csname#1tag\endcsname##1{\tagjj{##1}}\fi\fi\next}
\newif\ifshowstuffinmargin
\showstuffinmarginfalse
\def\jjtags{\ifx\shlhetal\relax 
  \else
\ifx\shlhetal\undefinedcontrolseq
\else
\showstuffinmargintrue
\ifx\all\relax\else\expandafter\t@gsjj\all,\fi\fi \fi
}

\def\tagjj#1{\realstag{#1}\oldmginpar{\zeigen{#1}}}
\def\citejj#1{\rechnen{#1}\mginpar{\zeigen{#1}}}     

\def\rechnen#1{\expandafter\ifx\csname stag#1\endcsname\relax ??\else
                           \csname stag#1\endcsname\fi}

\newdimen\theight

\def\marginfont{\sevenrm}

\def\trymarginbox#1{\setbox0=\hbox{\marginfont\hskip\marginshift #1}%
		\global\marginshift\wd0 
		\global\advance\marginshift\margindelta}

\def \oldmginpar#1{%
\ifvmode\setbox0\hbox to \hsize{\hfill\rlap{\marginfont\quad#1}}%
\ht0 0cm
\dp0 0cm
\box0\vskip-\baselineskip
\else 
             \vadjust{\trymarginbox{#1}%
		\ifdim\marginshift>\marginmax \global\marginshift\marginmin
			\trymarginbox{#1}%
                \fi
             \theight=\ht0
             \advance\theight by \dp0    \advance\theight by \lineskip
             \kern -\theight \vbox to \theight{\rightline{\rlap{\box0}}%
\vss}}\fi}

\newdimen\upordown
\global\upordown=8pt
\font\tinyfont=cmtt8 
\def\mginpar#1{\smash{\hbox to 0cm{\kern-10pt\raise7pt\hbox{\tinyfont #1}\hss}}}
\def\mginpar#1{{\hbox to 0cm{\kern-10pt\raise\upordown\hbox{\tinyfont #1}\hss}}\global\upordown-\upordown}


\def\t@gsoff#1,{\def\@{#1}\ifx\@\empty\let\next=\relax\else\let\next=\t@gsoff
   \def\@@{p}\ifx\@\@@\else
   \expandafter\gdef\csname#1cite\endcsname##1{\zeigen{##1}}
   \expandafter\gdef\csname#1page\endcsname##1{?}
   \expandafter\gdef\csname#1tag\endcsname##1{\zeigen{##1}}\fi\fi\next}
\def\verbatimtags{\showstuffinmarginfalse
\ifx\all\relax\else\expandafter\t@gsoff\all,\fi}
\def\zeigen#1{\hbox{$\scriptstyle\langle$}#1\hbox{$\scriptstyle\rangle$}}


\def\margintag#1{\ifshowstuffinmargin\oldmginpar{\zeigen{#1}}\fi}

\def\(#1){\edef\dot@g{\ifmmode\ifinner(\hbox{\noexpand\etag{#1}})
   \else\noexpand\eqno(\hbox{\noexpand\etag{#1}})\fi
   \else(\noexpand\ecite{#1})\fi}\dot@g}

\newif\ifbr@ck
\def\eat#1{}
\def\[#1]{\br@cktrue[\br@cket#1'X]}
\def\br@cket#1'#2X{\def\temp{#2}\ifx\temp\empty\let\next\eat
   \else\let\next\br@cket\fi
   \ifbr@ck\br@ckfalse\br@ck@t#1,X\else\br@cktrue#1\fi\next#2X}
\def\br@ck@t#1,#2X{\def\temp{#2}\ifx\temp\empty\let\neext\eat
   \else\let\neext\br@ck@t\def\temp{,}\fi
   \def\teemp{#1}\ifx\teemp\empty\else\rcite{#1}\fi\temp\neext#2X}
\def\resetbr@cket{\gdef\[##1]{[\rtag{##1}]}}
\def\references{\resetbr@cket\newsection References\par}

\newtoks\symb@ls\newtoks\s@mb@ls\newtoks\p@gelist\n@wcount\ftn@mber
    \ftn@mber=1\newif\ifftn@mbers\ftn@mbersfalse\newif\ifbyp@ge\byp@gefalse
\def\defm@rk{\ifftn@mbers\n@mberm@rk\else\symb@lm@rk\fi}
\def\n@mberm@rk{\xdef\m@rk{{\the\ftn@mber}}%
    \global\advance\ftn@mber by 1 }
\def\rot@te#1{\let\temp=#1\global#1=\expandafter\r@t@te\the\temp,X}
\def\r@t@te#1,#2X{{#2#1}\xdef\m@rk{{#1}}}
\def\b@@st#1{{$^{#1}$}}\def\str@p#1{#1}
\def\symb@lm@rk{\ifbyp@ge\rot@te\p@gelist\ifnum\expandafter\str@p\m@rk=1 
    \s@mb@ls=\symb@ls\fi\write\f@nsout{\number\count0}\fi \rot@te\s@mb@ls}
\def\byp@ge{\byp@getrue\n@wwrite\f@nsin\openin\f@nsin=\jobname.fns 
    \n@wcount\currentp@ge\currentp@ge=0\p@gelist={0}
    \re@dfns\closein\f@nsin\rot@te\p@gelist
    \n@wread\f@nsout\openout\f@nsout=\jobname.fns }
\def\m@kelist#1X#2{{#1,#2}}
\def\re@dfns{\ifeof\f@nsin\let\next=\relax\else\read\f@nsin to \f@nline
    \ifx\f@nline\v@idline\else\let\t@mplist=\p@gelist
    \ifnum\currentp@ge=\f@nline
    \global\p@gelist=\expandafter\m@kelist\the\t@mplistX0
    \else\currentp@ge=\f@nline
    \global\p@gelist=\expandafter\m@kelist\the\t@mplistX1\fi\fi
    \let\next=\re@dfns\fi\next}
\def\symbols#1{\symb@ls={#1}\s@mb@ls=\symb@ls} 
\def\bigsymbol{\textstyle}
\symbols{\bigsymbol\ast,\dagger,\ddagger,\sharp,\flat,\natural,\star}
\def\ftnumbers{\ftn@mberstrue} \def\ftsymbols{\ftn@mbersfalse}
\def\paginal{\byp@ge} \def\resetftnumbers{\ftn@mber=1}
\def\ftnote#1{\defm@rk\expandafter\expandafter\expandafter\footnote
    \expandafter\b@@st\m@rk{#1}}

\long\def\jump#1\endjump{}
\def\ssum{\mathop{\lower .1em\hbox{$\textstyle\Sigma$}}\nolimits}

\def\qed{\nobreak\kern 1em \vrule height .5em width .5em depth 0em}
\def\newneq{\hbox{\rlap{\hbox to 1\wd9{\hss$=$\hss}}\raise .1em 
   \hbox to 1\wd9{\hss$\scriptscriptstyle/$\hss}}}
\def\subsetne{\setbox9 = \hbox{$\subset$}\mathrel{\hbox{\rlap
   {\lower .4em \newneq}\raise .13em \hbox{$\subset$}}}}
\def\supsetne{\setbox9 = \hbox{$\subset$}\mathrel{\hbox{\rlap
   {\lower .4em \newneq}\raise .13em \hbox{$\supset$}}}}

\def\vbar{\mathchoice{\vrule height6.3ptdepth-.5ptwidth.8pt\kern-.8pt}
   {\vrule height6.3ptdepth-.5ptwidth.8pt\kern-.8pt}
   {\vrule height4.1ptdepth-.35ptwidth.6pt\kern-.6pt}
   {\vrule height3.1ptdepth-.25ptwidth.5pt\kern-.5pt}}
\def\f@dge{\mathchoice{}{}{\mkern.5mu}{\mkern.8mu}}
\def\b@c#1#2{{\rm \mkern#2mu\vbar\mkern-#2mu#1}}
\def\b@b#1{{\rm I\mkern-3.5mu #1}}
\def\b@a#1#2{{\rm #1\mkern-#2mu\f@dge #1}}
\def\bb#1{{\count4=`#1 \advance\count4by-64 \ifcase\count4\or\b@a A{11.5}\or
   \b@b B\or\b@c C{5}\or\b@b D\or\b@b E\or\b@b F \or\b@c G{5}\or\b@b H\or
   \b@b I\or\b@c J{3}\or\b@b K\or\b@b L \or\b@b M\or\b@b N\or\b@c O{5} \or
   \b@b P\or\b@c Q{5}\or\b@b R\or\b@a S{8}\or\b@a T{10.5}\or\b@c U{5}\or
   \b@a V{12}\or\b@a W{16.5}\or\b@a X{11}\or\b@a Y{11.7}\or\b@a Z{7.5}\fi}}

\catcode`\X=11 \catcode`\@=12




\let\thischap\jobname

\def\partof#1{\csname returnthe#1part\endcsname}
\def\chapof#1{\csname returnthe#1chap\endcsname}

\def\setchapter#1,#2,#3;{%
  \expandafter\def\csname returnthe#1part\endcsname{#2}%
  \expandafter\def\csname returnthe#1chap\endcsname{#3}%
}

\setchapter 300a,A,II.A;
\setchapter 300b,A,II.B;
\setchapter 300c,A,II.C;
\setchapter 300d,A,II.D;
\setchapter 300e,A,II.E;
\setchapter 300f,A,II.F;
\setchapter 300g,A,II.G;
\setchapter  E53,B,N;
\setchapter  88r,B,I;
\setchapter  600,B,III;
\setchapter  705,B,IV;
\setchapter  734,B,V;

\def\cprefix#1{
\edef\theotherpart{\partof{#1}}\edef\theotherchap{\chapof{#1}}%
\ifx\theotherpart\thispart
   \ifx\theotherchap\thischap 
    \else 
     \theotherchap%
    \fi
   \else 
     \theotherchap\fi}

\def\sectioncite[#1]#2{%
     \cprefix{#2}#1}

\edef\thispart{\partof{\thischap}}
\edef\thischap{\chapof{\thischap}}

\def\lastpage of '#1' is #2.{\expandafter\def\csname lastpage#1\endcsname{#2}}


\def\spuriousreset{}


\expandafter\ifx\csname citeadd.tex\endcsname\relax
\expandafter\gdef\csname citeadd.tex\endcsname{}
\else \message{Hey!  Apparently you were trying to
\string\input{citeadd.tex} twice.   This does not make sense.} 
\errmessage{Please edit your file (probably \jobname.tex) and remove
any duplicate ``\string\input'' lines}\endinput\fi

\sectno=-1   
\localtags
\jjtags
\NoBlackBoxes
\define\mr{\medskip\roster}
\define\sn{\smallskip\noindent}
\define\mn{\medskip\noindent}
\define\bn{\bigskip\noindent}
\define\ub{\underbar}

\define\ermn{\endroster\medskip\noindent}
\define\utt{\underset\tilde {}\to}

\define \nl{\newline}
\magnification=\magstep 1
\documentstyle{amsppt}

{    
\catcode`@11

\ifx\alicetwothousandloaded@\relax
  \endinput\else\global\let\alicetwothousandloaded@\relax\fi

\gdef\subjclass{\let\savedef@\subjclass
 \def\subjclass##1\endsubjclass{\let\subjclass\savedef@
   \toks@{\def\usualspace{{\rm\enspace}}\eightpoint}%
   \toks@@{##1\unskip.}%
   \edef\thesubjclass@{\the\toks@
     \frills@{{\noexpand\rm2000 {\noexpand\it Mathematics Subject
       Classification}.\noexpand\enspace}}%
     \the\toks@@}}%
  \nofrillscheck\subjclass}
} 


\expandafter\ifx\csname alice2jlem.tex\endcsname\relax
  \expandafter\xdef\csname alice2jlem.tex\endcsname{\the\catcode`@}
\else \message{Hey!  Apparently you were trying to
\string\input{alice2jlem.tex}  twice.   This does not make sense.}
\errmessage{Please edit your file (probably \jobname.tex) and remove
any duplicate ``\string\input'' lines}\endinput\fi

\expandafter\ifx\csname bib4plain.tex\endcsname\relax
  \expandafter\gdef\csname bib4plain.tex\endcsname{}
\else \message{Hey!  Apparently you were trying to \string\input
  bib4plain.tex twice.   This does not make sense.}
\errmessage{Please edit your file (probably \jobname.tex) and remove
any duplicate ``\string\input'' lines}\endinput\fi

\def\renewcommand{\newcommand}	       
\edef\cite{\the\catcode`@}%
\catcode`@ = 11
\let\@oldatcatcode = \cite
\chardef\@letter = 11
\chardef\@other = 12
%
%
%
%
\def\@innerdef#1#2{\edef#1{\expandafter\noexpand\csname #2\endcsname}}%
%
%
\@innerdef\@innernewcount{newcount}%
\@innerdef\@innernewdimen{newdimen}%
\@innerdef\@innernewif{newif}%
\@innerdef\@innernewwrite{newwrite}%
%
%
%
\def\@gobble#1{}%
%
%
%
\ifx\inputlineno\@undefined
   \let\@linenumber = \empty 
\else
   \def\@linenumber{\the\inputlineno:\space}%
\fi
%
%
%
\def\@futurenonspacelet#1{\def\cs{#1}%
   \afterassignment\@stepone\let\@nexttoken=
}%
\begingroup 
\def\\{\global\let\@stoken= }%
\\ 
\endgroup
\def\@stepone{\expandafter\futurelet\cs\@steptwo}%
\def\@steptwo{\expandafter\ifx\cs\@stoken\let\@@next=\@stepthree
   \else\let\@@next=\@nexttoken\fi \@@next}%
\def\@stepthree{\afterassignment\@stepone\let\@@next= }%
%
%
%
\def\@getoptionalarg#1{%
   \let\@optionaltemp = #1%
   \let\@optionalnext = \relax
   \@futurenonspacelet\@optionalnext\@bracketcheck
}%
%
%
\def\@bracketcheck{%
   \ifx [\@optionalnext
      \expandafter\@@getoptionalarg
   \else
      \let\@optionalarg = \empty
      \expandafter\@optionaltemp
   \fi
}%
\def\@@getoptionalarg[#1]{%
   \def\@optionalarg{#1}%
   \@optionaltemp
}%
%
%
%
\def\@nnil{\@nil}%
\def\@fornoop#1\@@#2#3{}%
\def\@for#1:=#2\do#3{%
   \edef\@fortmp{#2}%
   \ifx\@fortmp\empty \else
      \expandafter\@forloop#2,\@nil,\@nil\@@#1{#3}%
   \fi
}%
\def\@forloop#1,#2,#3\@@#4#5{\def#4{#1}\ifx #4\@nnil \else
       #5\def#4{#2}\ifx #4\@nnil \else#5\@iforloop #3\@@#4{#5}\fi\fi
}%
\def\@iforloop#1,#2\@@#3#4{\def#3{#1}\ifx #3\@nnil
       \let\@nextwhile=\@fornoop \else
      #4\relax\let\@nextwhile=\@iforloop\fi\@nextwhile#2\@@#3{#4}%
}%
%
%
%
\@innernewif\if@fileexists
\def\@testfileexistence{\@getoptionalarg\@finishtestfileexistence}%
\def\@finishtestfileexistence#1{%
   \begingroup
      \def\extension{#1}%
      \immediate\openin0 =
         \ifx\@optionalarg\empty\jobname\else\@optionalarg\fi
         \ifx\extension\empty \else .#1\fi
         \space
      \ifeof 0
         \global\@fileexistsfalse
      \else
         \global\@fileexiststrue
      \fi
      \immediate\closein0
   \endgroup
}%
%
%
%
%
\def\bibliographystyle#1{%
   \@readauxfile
   \@writeaux{\string\bibstyle{#1}}%
}%
\let\bibstyle = \@gobble
%
%
\let\bblfilebasename = \jobname
\def\bibliography#1{%
   \@readauxfile
   \@writeaux{\string\bibdata{#1}}%
   \@testfileexistence[\bblfilebasename]{bbl}%
   \if@fileexists
      \nobreak
      \@readbblfile
   \fi
}%
\let\bibdata = \@gobble
%
%
\def\nocite#1{%
   \@readauxfile
   \@writeaux{\string\citation{#1}}%
}%
\@innernewif\if@notfirstcitation
%
%
\def\cite{\@getoptionalarg\@cite}%
%
%
\def\@cite#1{%
   \let\@citenotetext = \@optionalarg
   \printcitestart
   \nocite{#1}%
   \@notfirstcitationfalse
   \@for \@citation :=#1\do
   {%
      \expandafter\@onecitation\@citation\@@
   }%
   \ifx\empty\@citenotetext\else
      \printcitenote{\@citenotetext}%
   \fi
   \printcitefinish
}%
\newif\ifweareinprivate
\weareinprivatetrue
\ifx\shlhetal\undefinedcontrolseq\weareinprivatefalse\fi
\ifx\shlhetal\relax\weareinprivatefalse\fi
\def\@onecitation#1\@@{%
   \if@notfirstcitation
      \printbetweencitations
   \fi
   \expandafter \ifx \csname\@citelabel{#1}\endcsname \relax
      \if@citewarning
         \message{\@linenumber Undefined citation `#1'.}%
      \fi
     \ifweareinprivate
      \expandafter\gdef\csname\@citelabel{#1}\endcsname{%
\strut 
\vadjust{\vskip-\dp\strutbox
\vbox to 0pt{\vss\parindent0cm \leftskip=\hsize 
\advance\leftskip3mm
\advance\hsize 4cm\strut\openup-4pt 
\rightskip 0cm plus 1cm minus 0.5cm ?  #1 ?\strut}}
         {\tt
            \escapechar = -1
            \nobreak\hskip0pt\pfeilsw
            \expandafter\string\csname#1\endcsname
             \pfeilso
            \nobreak\hskip0pt
         }%
      }%
     \else  
      \expandafter\gdef\csname\@citelabel{#1}\endcsname{%
            {\tt\expandafter\string\csname#1\endcsname}
      }%
     \fi  
   \fi
   \csname\@citelabel{#1}\endcsname
   \@notfirstcitationtrue
}%
%
%
\def\@citelabel#1{b@#1}%
%
%
\def\@citedef#1#2{\expandafter\gdef\csname\@citelabel{#1}\endcsname{#2}}%
%
%
%
\def\@readbblfile{%
   \ifx\@itemnum\@undefined
      \@innernewcount\@itemnum
   \fi
   \begingroup
      \def\begin##1##2{%
         \setbox0 = \hbox{\biblabelcontents{##2}}%
         \biblabelwidth = \wd0
      }%
      \def\end##1{}
      %
      %
      \@itemnum = 0
      \def\bibitem{\@getoptionalarg\@bibitem}%
      \def\@bibitem{%
         \ifx\@optionalarg\empty
            \expandafter\@numberedbibitem
         \else
            \expandafter\@alphabibitem
         \fi
      }%
      \def\@alphabibitem##1{%
         \expandafter \xdef\csname\@citelabel{##1}\endcsname {\@optionalarg}%
         \ifx\biblabelprecontents\@undefined
            \let\biblabelprecontents = \relax
         \fi
         \ifx\biblabelpostcontents\@undefined
            \let\biblabelpostcontents = \hss
         \fi
         \@finishbibitem{##1}%
      }%
      \def\@numberedbibitem##1{%
         \advance\@itemnum by 1
         \expandafter \xdef\csname\@citelabel{##1}\endcsname{\number\@itemnum}%
         \ifx\biblabelprecontents\@undefined
            \let\biblabelprecontents = \hss
         \fi
         \ifx\biblabelpostcontents\@undefined
            \let\biblabelpostcontents = \relax
         \fi
         \@finishbibitem{##1}%
      }%
      \def\@finishbibitem##1{%
         \biblabelprint{\csname\@citelabel{##1}\endcsname}%
         \@writeaux{\string\@citedef{##1}{\csname\@citelabel{##1}\endcsname}}%
         \ignorespaces
      }%
      %
      %
      \let\em = \bblem
      \let\newblock = \bblnewblock
      \let\sc = \bblsc
      \frenchspacing
      \clubpenalty = 4000 \widowpenalty = 4000
      \tolerance = 10000 \hfuzz = .5pt
      \everypar = {\hangindent = \biblabelwidth
                      \advance\hangindent by \biblabelextraspace}%
      \bblrm
      \parskip = 1.5ex plus .5ex minus .5ex
      \biblabelextraspace = .5em
      \bblhook
      \input \bblfilebasename.bbl
   \endgroup
}%
%
%
\@innernewdimen\biblabelwidth
\@innernewdimen\biblabelextraspace
%
%
%
\def\biblabelprint#1{%
   \noindent
   \hbox to \biblabelwidth{%
      \biblabelprecontents
      \biblabelcontents{#1}%
      \biblabelpostcontents
   }%
   \kern\biblabelextraspace
}%
%
%
%
\def\biblabelcontents#1{{\bblrm [#1]}}%
%
%
\def\bblrm{\rm}%
%
%
\def\bblem{\it}%
%
%
\def\bblsc{\ifx\@scfont\@undefined
              \font\@scfont = cmcsc10
           \fi
           \@scfont
}%
%
%
\def\bblnewblock{\hskip .11em plus .33em minus .07em }%
%
%
\let\bblhook = \empty
%
%
%
\def\printcitestart{[}
\def\printcitefinish{]}
\def\printbetweencitations{, }
\def\printcitenote#1{, #1}
%
%
%
\let\citation = \@gobble
%
%
%
\@innernewcount\@numparams
%
%
\def\newcommand#1{%
   \def\@commandname{#1}%
   \@getoptionalarg\@continuenewcommand
}%
%
%
\def\@continuenewcommand{%
   \@numparams = \ifx\@optionalarg\empty 0\else\@optionalarg \fi \relax
   \@newcommand
}%
%
%
\def\@newcommand#1{%
   \def\@startdef{\expandafter\edef\@commandname}%
   \ifnum\@numparams=0
      \let\@paramdef = \empty
   \else
      \ifnum\@numparams>9
         \errmessage{\the\@numparams\space is too many parameters}%
      \else
         \ifnum\@numparams<0
            \errmessage{\the\@numparams\space is too few parameters}%
         \else
            \edef\@paramdef{%
               \ifcase\@numparams
                  \empty  No arguments.
               \or ####1%
               \or ####1####2%
               \or ####1####2####3%
               \or ####1####2####3####4%
               \or ####1####2####3####4####5%
               \or ####1####2####3####4####5####6%
               \or ####1####2####3####4####5####6####7%
               \or ####1####2####3####4####5####6####7####8%
               \or ####1####2####3####4####5####6####7####8####9%
               \fi
            }%
         \fi
      \fi
   \fi
   \expandafter\@startdef\@paramdef{#1}%
}%
%
%
%
%
\def\@readauxfile{%
   \if@auxfiledone \else 
      \global\@auxfiledonetrue
      \@testfileexistence{aux}%
      \if@fileexists
         \begingroup
            \endlinechar = -1
            \catcode`@ = 11
            \input \jobname.aux
         \endgroup
      \else
         \message{\@undefinedmessage}%
         \global\@citewarningfalse
      \fi
      \immediate\openout\@auxfile = \jobname.aux
   \fi
}%
%
%
\newif\if@auxfiledone
\ifx\noauxfile\@undefined \else \@auxfiledonetrue\fi
%
%
%
%
\@innernewwrite\@auxfile
\def\@writeaux#1{\ifx\noauxfile\@undefined \write\@auxfile{#1}\fi}%
%
%
%
\ifx\@undefinedmessage\@undefined
   \def\@undefinedmessage{No .aux file; I won't give you warnings about
                          undefined citations.}%
\fi
%
%
\@innernewif\if@citewarning
\ifx\noauxfile\@undefined \@citewarningtrue\fi
%
%
%
\catcode`@ = \@oldatcatcode

\def\pfeilso{\leavevmode
            \vrule width 1pt height9pt depth 0pt\relax
           \vrule width 1pt height8.7pt depth 0pt\relax
           \vrule width 1pt height8.3pt depth 0pt\relax
           \vrule width 1pt height8.0pt depth 0pt\relax
           \vrule width 1pt height7.7pt depth 0pt\relax
            \vrule width 1pt height7.3pt depth 0pt\relax
            \vrule width 1pt height7.0pt depth 0pt\relax
            \vrule width 1pt height6.7pt depth 0pt\relax
            \vrule width 1pt height6.3pt depth 0pt\relax
            \vrule width 1pt height6.0pt depth 0pt\relax
            \vrule width 1pt height5.7pt depth 0pt\relax
            \vrule width 1pt height5.3pt depth 0pt\relax
            \vrule width 1pt height5.0pt depth 0pt\relax
            \vrule width 1pt height4.7pt depth 0pt\relax
            \vrule width 1pt height4.3pt depth 0pt\relax
            \vrule width 1pt height4.0pt depth 0pt\relax
            \vrule width 1pt height3.7pt depth 0pt\relax
            \vrule width 1pt height3.3pt depth 0pt\relax
            \vrule width 1pt height3.0pt depth 0pt\relax
            \vrule width 1pt height2.7pt depth 0pt\relax
            \vrule width 1pt height2.3pt depth 0pt\relax
            \vrule width 1pt height2.0pt depth 0pt\relax
            \vrule width 1pt height1.7pt depth 0pt\relax
            \vrule width 1pt height1.3pt depth 0pt\relax
            \vrule width 1pt height1.0pt depth 0pt\relax
            \vrule width 1pt height0.7pt depth 0pt\relax
            \vrule width 1pt height0.3pt depth 0pt\relax}

\def\pfeilsw{ \leavevmode 
            \vrule width 1pt height0.3pt depth 0pt\relax
            \vrule width 1pt height0.7pt depth 0pt\relax
            \vrule width 1pt height1.0pt depth 0pt\relax
            \vrule width 1pt height1.3pt depth 0pt\relax
            \vrule width 1pt height1.7pt depth 0pt\relax
            \vrule width 1pt height2.0pt depth 0pt\relax
            \vrule width 1pt height2.3pt depth 0pt\relax
            \vrule width 1pt height2.7pt depth 0pt\relax
            \vrule width 1pt height3.0pt depth 0pt\relax
            \vrule width 1pt height3.3pt depth 0pt\relax
            \vrule width 1pt height3.7pt depth 0pt\relax
            \vrule width 1pt height4.0pt depth 0pt\relax
            \vrule width 1pt height4.3pt depth 0pt\relax
            \vrule width 1pt height4.7pt depth 0pt\relax
            \vrule width 1pt height5.0pt depth 0pt\relax
            \vrule width 1pt height5.3pt depth 0pt\relax
            \vrule width 1pt height5.7pt depth 0pt\relax
            \vrule width 1pt height6.0pt depth 0pt\relax
            \vrule width 1pt height6.3pt depth 0pt\relax
            \vrule width 1pt height6.7pt depth 0pt\relax
            \vrule width 1pt height7.0pt depth 0pt\relax
            \vrule width 1pt height7.3pt depth 0pt\relax
            \vrule width 1pt height7.7pt depth 0pt\relax
            \vrule width 1pt height8.0pt depth 0pt\relax
            \vrule width 1pt height8.3pt depth 0pt\relax
            \vrule width 1pt height8.7pt depth 0pt\relax
            \vrule width 1pt height9pt depth 0pt\relax
      }


\def\widestnumber#1#2{}

\def\citewarning#1{\ifx\shlhetal\relax 
    \else
    \par{#1}\par
    \fi
}

\def\rm{\fam0 \tenrm}

\def\fakesubhead#1\endsubhead{\bigskip\noindent{\bf#1}\par}



%
%
%

%

\font\textrsfs=rsfs10
\font\scriptrsfs=rsfs7
\font\scriptscriptrsfs=rsfs5

\newfam\rsfsfam
\textfont\rsfsfam=\textrsfs
\scriptfont\rsfsfam=\scriptrsfs
\scriptscriptfont\rsfsfam=\scriptscriptrsfs

\edef\oldcatcodeofat{\the\catcode`\@}
\catcode`\@11

\def\Cal@@#1{\noaccents@ \fam \rsfsfam #1}

\catcode`\@\oldcatcodeofat


\expandafter\ifx \csname margininit\endcsname \relax\else\margininit\fi

\long\def\red#1\endred{}
\long\def\green#1\endgreen{}
\long\def\blue#1\endblue{}
\long\def\private#1\endprivate{}

\def\endred{ \unmatched endred! }
\def\endgreen{ \unmatched endgreen! }
\def\endblue{ \unmatched endblue! }
\def\endprivate{ \unmatched endprivate! }

\ifx\latexcolors\undefinedcs\def\latexcolors{}\fi

\def\emptycs{}
\def\evaluatelatexcolors{%
        \ifx\latexcolors\emptycs\else
        \expandafter\xxevaluate\latexcolors\xxfertig\evaluatelatexcolors\fi}
\def\xxevaluate#1,#2\xxfertig{\setupthiscolor{#1}%
        \def\latexcolors{#2}}


\font\smallfont=cmsl7
\def\rutgerscolor{\ifmmode\else\endgraf\fi\smallfont
\advance\leftskip0.5cm\relax}
\def\setupthiscolor#1{\edef\tmptmpcs{\noexpand\bgroup\noexpand\rutgerscolor
\noexpand\def\noexpand\currentcolor{#1}%
\noexpand}%
\expandafter\let\csname#1\endcsname\tmptmpcs
\def\tmptmpcs{\checkColorUnmatched{#1}\popthecolor}
\expandafter\let\csname end#1\endcsname\tmptmpcs}

\def\checkColorUnmatched#1{\def\expectcolor{#1}%
    \ifx\expectcolor\currentcolor   
    \else \edef\failhere{\noexpand\tryingToClose '\currentcolor' with end\expectcolor}\failhere\fi}

\def\currentcolor{???}

\def\popthecolor{\ifmmode\else\endgraf\fi\egroup}

\expandafter\def\csname#1\endcsname{}

\evaluatelatexcolors

 \let\outerhead\head
 \def\head{\innerhead}
 \let\innerhead\outerhead

 \let\outersubhead\subhead
 \def\subhead{\innersubhead}
 \let\innersubhead\outersubhead

 \let\outersubsubhead\subsubhead
 \def\subsubhead{\innersubsubhead}
 \let\innersubsubhead\outersubsubhead

 \let\outerproclaim\proclaim
 \def\proclaim{\innerproclaim}
 \let\innerproclaim\outerproclaim

 %
 %
 %
 %

\def\demo#1{\medskip\noindent{\it #1.\/}}
\def\enddemo{\smallskip}

\def\remark#1{\medskip\noindent{\it #1.\/}}
\def\endremark{\smallskip}

\pageheight{8.5truein}
\topmatter
\title {The spectrum of characters of ultrafilters on $\omega$ } \endtitle
\author {Saharon Shelah \thanks {\null\newline I would like to thank 
Alice Leonhardt for the beautiful typing. \null\newline
 Partially supported by the Binational Science Foundation and the
Canadian Research Chair; 613-943-9382. Publication 846} 
\endthanks} \endauthor 

\affil{The Hebrew University of Jerusalem \\
Einstein Institute of Mathematics \\
Edmond J. Safra Campus, Givat Ram \\
Jerusalem 91904, Israel
 \medskip
 Department of Mathematics \\
 Hill Center-Busch Campus \\
  Rutgers, The State University of New Jersey \\
 110 Frelinghuysen Road \\
 Piscataway, NJ 08854-8019 USA} \endaffil


\abstract  We show the consistency of statement: ``the set of regular cardinals
which are the character of some ultrafilter on $\omega$ is not convex".
We also deal with the set of $\pi$-characters of ultrafilters 
on $\omega$. \endabstract
\endtopmatter
\document

\head {\S0 Introduction} \endhead  \resetall \sectno=0
 \spuriousreset
\bigskip

Some cardinal invariants of the continuum are actually the minimum
of a natural set of cardinals $\le 2^{\aleph_0}$ which can be called
the spectrum of the invariant.  Such a
case is Sp$_\chi$, the set of characters $\chi(D)$ of non-principal
ultrafilters on $\omega$ (the minimal number of generators).  On the
history see \cite{BnSh:642}; there this spectrum and others were 
investigated and it was asked if Sp$_\chi$ can be non-convex (formally
\scite{642.1A}(2) below).

The main result is \scite{642.2}, it solves the problem (starting with a
measurable).
This was presented in a conference in
honor of Juhasz, quite fitting as he had started the investigation of
consistency on $\chi(D)$.  In \S2 we note what we can say on the
strict $\pi$-character of ultrafilters.

The investigation is continued in \cite{Sh:915} trying to get more
``disorderly" behaviours in smaller cardinals and in particular
answering negatively the original question, \scite{642.1A}(2).
\bn
Recall
\definition{\stag{642.1} Definition}  1) Sp$_\chi = \text{ Sp}(\chi)$ 
is the set of cardinals $\theta$ such that: $\theta = \chi(D)$ for 
some non-principal ultrafilter $D$ on $\omega$ where \nl
2) For $D$ an ultrafilter on $\omega$ let $\theta = \chi(D)$ be the
minimal cardinality $\theta$ such that $D$ is generated by 
some family of $\theta$ members, i.e. Min$\{|{\Cal A}|:{\Cal A}
\subseteq D$ and $(\forall B \in D)(\exists A \in {\Cal A})[A
\subseteq^* B]\}$, it does not matter if we use ``$A \subseteq B$". 
\enddefinition
\bn
Now, Brendle and Shelah \cite[Problem 5]{BnSh:642}, asked the question
formulated in \scite{642.1A}(2) below, but it seems to me now 
that the question is really \scite{642.1A}(1)+(3). 
\bn
\margintag{642.1A}\ub{\stag{642.1A} Problem}  1) Can Sp$(\chi) \cap \text{ Reg}$ have gaps,
i.e., can it be that $\theta < \mu < \lambda$ are regular, $\theta \in
\text{ Sp}(\chi),\mu \notin \text{ Sp}(\chi),\lambda \in \text{
Sp}(\chi)$? 
\nl
2) In particular does $\aleph_1,\aleph_3 \in \text{ Sp}(\chi)$ imply
$\aleph_2 \in \text{ Sp}(\chi)$? 
\nl
3) Are there any restrictions on Sp$(\chi) \cap \text{ Reg}$?

We thank the referee for helpful comments and in particular
\scite{642.2A}(1).
\bn
\ub{Discussion}:  This rely on \cite[\S4]{Sh:700}, there is no point to
repeat it but we try to give a description.

Let $S = \{\alpha < \lambda:\text{cf}(\alpha) \ne \kappa\}$ or any unbounded
subset of it.  We define (\cite[4.3]{Sh:700}) the class ${\frak K} =
{\frak K}_{\lambda,S}$ of objects ${\frak t}$ approximating our final
forcing.  Each ${\frak t} \in K$ consists mainly of a finite support
iteration $\langle \Bbb P^{\frak t}_i,{\utt {\Bbb Q}^{\frak t}_i}:i <
\mu\rangle$ of c.c.c. forcing of cardinality $\le \lambda$ with limit
$\Bbb P^*_{\frak t} = \Bbb P^{\frak t} = \Bbb P^{\frak t}_\mu$, but also
${\underset\tilde {}\to \tau^{\frak t}_i} \, (i < \mu)$ of $\Bbb
Q^{\frak t}_i$ satisfying a strong version of the c.c.c. and for $i
\in S$, also ${\underset\tilde {}\to D^{\frak t}_i}$, a $\Bbb
P^t_i$-name of a non-principal ultrafilter on $\omega$ from which
$\Bbb Q^{\frak t}_i$ is nicely defined and ${\underset\tilde {}\to
A^{\frak t}_i}$, a ${\utt {\Bbb Q}^{\frak t}_i}$-name (so $\Bbb
P^t_{i+1}$-name) of a pseudo-intersection (and $\Bbb Q_i,i \in S$,
nicely defined) of ${\underset\tilde {}\to D^t_i}$ such that $i < j
\in S \Rightarrow {\underset\tilde {}\to A^{\frak t}_i} \in 
D^{\frak t}_j$. So $\{{\underset\tilde {}\to A_i}:i \in S\}$ witness ${\frak u}
\le \mu$ in $\bold V^{\Bbb P_{\frak t}}$; not necessarily we have to
use non-nicely defined $\Bbb Q_i$, though for $i \in S$ we do.

The order $\le_{\frak K}$ is natural order, we prove the existence of
the so-called canonical limit.

Now a major point of \cite{Sh:700} is: for ${\frak s} \in {\frak K}$
letting ${\Cal D}$ be a uniform $\kappa$-complete ultrafilter on
$\kappa$, (or just $\kappa_1$-complete $\aleph_0 < \theta < \kappa$),
we can consider ${\frak t} = {\frak s}^\kappa/{\Cal D}$; by \L os
theorem, more exactly by Hanf's Ph.D. Thesis, 
(the parallel of) \L os theorem for $\Bbb
L_{\kappa_1,\kappa_2}$ apply, it gives that ${\frak t} \in {\frak K}$,
well if $\lambda = \lambda^\kappa/{\Cal D}$; and moreover ${\frak s}
\le_{\frak K} {\frak t}$ under the canonical embedding.

The effect is that, e.g. being ``a linear order having cofinality $\theta \ne
\kappa$" is preserved, even by the same witness whereas having
cardinality $\theta < \lambda$ is not and sets of cardinality $\ge \kappa$ are
increased.  As ${\frak d}$ is the cofinality (not of a linear order
but) of a partial order there are complications,
anyhow as ${\frak d}$ is defined by cofinality whereas ${\frak a}$ by
cardinality of sets this helps in \cite{Sh:700}, noting that as we
deal with c.c.c. forcing, reals are represented by $\omega$-sequences
of conditions, the relevant thing are preserved.  So we use a
$\le_{\frak K}$-increasing sequence $\langle {\frak t}_\alpha:\alpha
\le \lambda\rangle$ such that for unboundedly many $\alpha <
\lambda,{\frak t}_{\alpha +1}$ is essentially $({\frak
t}^\alpha_\alpha)^\kappa/{\Cal D}$.

What does ``nice" $\Bbb Q = \Bbb Q(D)$, for $D$ a non-principal
ultrafilter over $\omega$ mean?  We need that
\mr
\item "{$(\alpha)$}"  $\Bbb Q$ satisfies a strong version of the
c.c.c.
\sn
\item "{$(\beta)$}"  the definition commute with the ultra-power used
\sn
\item "{$(\gamma)$}"  if $\Bbb P$ is a forcing notion then we can
extend $D$ to an ultrafilter ${\underset\tilde {}\to D^+}$ for every
(or at least some) $\Bbb P$-name of an ultrafilter ${\underset\tilde
{}\to D}$ extending $D$ we have $\Bbb Q(D) \lessdot \Bbb P * \Bbb
Q({\underset\tilde {}\to D^+})$ (used for the existence of canonical
limit).
\ermn
Such a forcing is combining Laver forcing and Mathias forcing for an
ultrafilter $D$ on $\omega$, that is:  if $p \in D$ iff $p$ is a subtree of
$\omega$ with trunk tr$(p) \in p$ such that for $\eta \in p$ we have
$\ell g(\eta) < \ell g(\text{tr}(p)) \rightarrow (\exists!n)(\eta
\char 94 \langle n\rangle \in p)$ and $\ell g(\eta) \ge \ell
g(\text{tr}(p)) \Rightarrow \{n:\eta \char 94 \langle n\rangle \in p\}
\in D$. 
\bigskip

\head {\S1 Using measurables and FS iterations with non-transitive
memory} \endhead  \resetall 
 \spuriousreset
\bigskip

We use \cite{Sh:700} in \scite{642.2} heavily.  
We use measurables (we could have used extenders
to get more).  The question on $\aleph_1,\aleph_2,\aleph_3$, i.e.
Problem \scite{642.1A}(2) remains open.
\bigskip

\proclaim{\stag{642.2} Theorem}   There is a 
{\rm c.c.c.} forcing notion $\Bbb
P$ of cardinality $\lambda$ such that in $\bold V^{\Bbb P}$ we have 
${\frak a} =
\lambda,{\frak b} = {\frak d} = \mu,{\frak u} = \mu,\{\mu,\lambda\}
\subseteq { \text{\rm Sp\/}}_\chi$ but 
$\kappa_2 \notin { \text{\rm Sp\/}}(\chi)$ \ub{if}
\mr
\item "{$\circledast_1$}"  $\kappa_1,\kappa_2$ are measurable and 
$\kappa_1 < \mu =
{ \text{\rm cf\/}}(\mu) < \kappa_2 < \lambda = \lambda^\mu =
\lambda^{\kappa_2} = { \text{\rm cf\/}}(\lambda)$.
\endroster
\endproclaim
\bigskip

\demo{Proof}  Let ${\Cal D}_\ell$ be a normal ultrafilter on
$\kappa_\ell$ for $\ell=1,2$.   Repeat \cite[\S4]{Sh:700} with
$(\kappa_1,\mu,\lambda)$ here standing for $(\kappa,\mu,\lambda)$
there, getting ${\frak t}_\alpha \in {\frak K}$ for $\alpha \le
\lambda$ which is
$\le_{\frak K}$-increasing and letting $\Bbb P^\alpha_i = 
\Bbb P^{{\frak t}_\alpha}_i$ we have $\bar{\Bbb Q}^\alpha = 
\langle \Bbb P^\alpha_\varepsilon:\varepsilon < \mu\rangle$ is
a $\lessdot$-increasing continuous sequence of c.c.c. 
forcing notions, $\Bbb P^\alpha_\mu = \Bbb P^\alpha = 
\Bbb P_{{\frak t}_\alpha} := 
\text{ Lim}(\bar{\Bbb Q}^\alpha) = \cup\{\Bbb P^\alpha_\varepsilon:
\varepsilon < \mu\}$ but add the demand
that for unboundedly many $\alpha < \lambda$
\mr
\item "{$\boxtimes^1_\alpha$}"  $\Bbb P^{\alpha +1}$ is isomorphic to the
ultrapower $(\Bbb P^\alpha)^{\kappa_2}/{\Cal D}_2$, by an isomorphism extending
the canonical embedding.
\ermn
More explicitly we choose ${\frak t}_\alpha$ by induction on $\alpha \le
\lambda$ such that
\mr
\item "{$\circledast_1$}"  $(a) \quad {\frak t}_\alpha \in {\frak K}$,
see Definition \cite[4.3]{Sh:700} so the forcing notion $\Bbb P^{{\frak
t}_\alpha}_i$ for $i \le \mu$
\nl

\hskip25pt is well defined and is $\lessdot$-increasing with $i$
\sn
\item "{${{}}$}"  $(b) \quad \langle {\frak t}_\beta:\beta \le \alpha
\rangle$ is $\le_{\frak K}$-increasing continuous which means that:
\sn
\item "{${{}}$}"  \hskip20pt $(\alpha) \quad \gamma \le \beta \le
\alpha \Rightarrow {\frak t}_\gamma \le_{\frak K} {\frak t}_\beta$, see
Definition \cite[4.6]{Sh:700}(1) so $\Bbb P^{{\frak t}_\gamma}_i
 \lessdot \Bbb P^{{\frak t}_\beta}_i$
\nl

\hskip25pt  for $i \le \mu$
 \sn
\item "{${{}}$}"  \hskip20pt $(\beta) \quad$ if $\alpha$ is a limit
ordinal then ${\frak t}_\alpha$ is a canonical $\le_{\frak K}$-u.b.
\nl

\hskip25pt  of $\langle {\frak t}_\beta:\beta < \alpha \rangle$, 
\nl

\hskip25pt see Definition \cite[4.6]{Sh:700}(2)
\sn
\item "{${{}}$}"  $(c) \quad$ if $\alpha = \beta +1$ and cf$(\beta) \ne
\kappa_2$ \ub{then} ${\frak t}_\alpha$ is essentially 
${\frak t}^{\kappa_1}_\beta/{\Cal D}_1$
\nl

\hskip25pt (i.e. we have to identify $\Bbb P^{{\frak
t}_\beta}_\varepsilon$ with its image under the 
\nl

\hskip25pt canonical embedding of
it into $(\Bbb P^{{\frak t}_\beta}_\varepsilon)^{\kappa_1}/{\Cal D}_1$, in 
\nl

\hskip25pt particular this holds for $\varepsilon =\mu$, see
Subclaim \cite[4.9]{Sh:700})
\sn
\item "{${{}}$}"  $(d) \quad$ if $\alpha= \beta +1$ and
cf$(\beta)=\kappa_2$ \ub{then} ${\frak t}_\alpha$ is essentially ${\frak
t}^{\kappa_2}_\beta/{\Cal D}_2$.
\nl
So we need
\sn
\item "{$\circledast_2$}"  Subclaim \cite[4.9]{Sh:700} applies also to
the ultrapower ${\frak t}^{\kappa_2}_\beta/D$.
\nl
[Why?  The same proof applies as $\mu^{\kappa_2}/{\Cal D}_2 = \mu$, i.e., the
canonical embedding of $\mu$ into $\mu^{\kappa_2}/{\Cal D}_2$ is one-to-one
and onto (and $\lambda^{\kappa_1}/{\Cal D}_1 = 
\lambda^{\kappa_2}/{\Cal D}_2 = \lambda$, of course).]
\ermn
Let $\Bbb P^\alpha_\varepsilon = \Bbb P^{{\frak
t}_\alpha}_\varepsilon$ for $\varepsilon \le \mu$ so 
$\Bbb P^\alpha = \cup\{\Bbb P^\alpha_\varepsilon:\varepsilon <
\mu\}$ and $\Bbb P = \Bbb P^\lambda$.  It is proved in \cite[4.10]{Sh:700} 
that in $\bold V^{\Bbb P}$, by the 
construction, $\mu \in \text{ Sp}(\chi),{\frak a} \le \lambda$ and
${\frak u} = \mu,2^{\aleph_0} = \lambda$.  
By \cite[4.11]{Sh:700} we have ${\frak a} \ge \lambda$ hence
${\frak a} = \lambda$, and always $2^{\aleph_0} \in \text{ Sp}(\chi)$
hence $\lambda = 2^{\aleph_0} \in \text{ Sp}(\chi)$.  So what is left to be
proved is $\kappa_2 \notin \text{ Sp}(\chi)$.  Assume
toward contradiction that $p^* \Vdash ``\underset\tilde {}\to D$
is a non-principal ultrafilter on $\omega$ and 
$\chi (\underset\tilde {}\to D) = \kappa_2$ and let it be exemplified 
by $\langle {\underset\tilde {}\to  A_\varepsilon}:
\varepsilon < \kappa_2 \rangle"$.
\sn
Without loss of generality 
$p^* \Vdash_{\Bbb P} ``{\underset\tilde {}\to
A_\varepsilon} \in \underset\tilde {}\to D$ does not belong to the
filter on $\omega$ generated by 
$\{{\underset\tilde {}\to A_\zeta}:\zeta < \varepsilon\} \cup \{\omega
\backslash n:n < \omega\}$, 
for each $\varepsilon < \kappa_2$ and trivially also
$\omega \backslash {\underset\tilde {}\to A_\varepsilon}$ does not
belong to this filter".

As $\lambda$ is regular $> \kappa_2$ and 
the forcing notion $\Bbb P^\lambda$ satisfies the
c.c.c., clearly for some $\alpha < \lambda$ we have
$p^* \in \Bbb P^\alpha$ and $\varepsilon < \kappa_2
\Rightarrow {\underset\tilde {}\to A_\varepsilon}$ is a
$\Bbb P^\alpha$-name.\nl
So for every $\beta \in [\alpha,\lambda)$ we have
\mr
\item "{$\boxtimes^2_\beta$}"  
$p^* \Vdash_{{\Bbb P}^\beta}$ ``for each $i < \kappa_2$ the set
${\underset\tilde {}\to A_i} \in  [\omega]^{\aleph_0}$ is 
not in the filter on $\omega$ which
$\{{\underset\tilde {}\to A_j}:j < i\} \cup\{\omega \backslash n:n <
\omega\}$ generates, and also the complement 
of ${\underset\tilde {}\to A_i}$ is
not in this filter (as $\underset\tilde {}\to D$ exemplifies this)".
\ermn
But for some such $\beta$, the statement $\boxtimes^1_\beta$ holds,
i.e. $\circledast_1(d)$ apply, 
so in $\Bbb P^{\beta+1}$ which essentially is a
$(\Bbb P^\beta)^{\kappa_2}/{\Cal D}_2$ we get a contradiction.
That is, let $\bold j_\beta$ be an isomorphism from $\Bbb P^{\beta
+1}$ onto $(\Bbb P^\beta)^{\kappa_2}/{\Cal D}_2$ which extends the canonical
embedding of $\Bbb P^\beta$ into $(\Bbb P^\beta)^{\kappa_2}/{\Cal D}_2$.
Now $\bold j_\beta$ induces a map $\hat{\bold j}_\beta$ from the set of
$\Bbb P^{\beta +1}$-names of subsets of $\omega$ into the set of $(\Bbb
P^\beta)^{\kappa_2}/{\Cal D}_2$-names of subsets of $\omega$, and
 let ${\underset\tilde {}\to A^*} = \hat{\bold j}^{-1}_\beta
(\langle {\underset\tilde {}\to A_i}:i < \kappa_2
\rangle/{\Cal D}_2)$ so $p^* \Vdash_{{\Bbb P}^{\beta +1}} ``{\underset\tilde
{}\to A^*} \in [\omega]^{\aleph_0}$ and the sets
${\underset\tilde {}\to A^*},\omega \backslash {\underset\tilde {}\to
A^*}$ do not include any finite intersection of
$\{{\underset\tilde {}\to A_\varepsilon}:\varepsilon < \kappa_2\} \cup
\{\omega \backslash n:n < \omega\}"$.  So $p^* \Vdash_{{\Bbb P}^{\beta +1}}
``\{{\underset\tilde {}\to A_\varepsilon}:\varepsilon < \kappa_2\}$ does
not generate an ultrafilter on $\omega$" but $\Bbb P^{\beta +1} \lessdot
\Bbb P$, contradiction.  \hfill$\square_{\scite{642.2}}$
\enddemo
\bigskip

\remark{\stag{642.2D} Remark}  1) As the referree pointed out we can in
\scite{642.2}, if we waive ``${\frak u} < {\frak a}$" we can forget
$\kappa_1$ (and ${\Cal D}_1$) so not taking ultra-powers by ${\Cal
D}_1$, so $\mu = \aleph_0$ is allowed, but we have to start with
${\frak t}_0$ such that $\Bbb P^{{\frak t}_0}_0$ is adding $\kappa_2$-Cohen.
\nl
2) Moreover, in this case we can demand that ${\underset\tilde {}\to
{\Bbb Q}^{\frak t}_\alpha} = 
{\underset\tilde {}\to {\Bbb Q}}
({\underset\tilde {}\to D^{\frak t}_\alpha})$ and so we do not
need the ${\underset\tilde {}\to t^{\frak t}_\alpha}$.  
Still this way was taken in \cite[\S1]{Sh:915}.
But this gain in simplicity has a price in lack of flexibility in choosing
 the ${\frak t}$.  
 We use this mildly in \S2; mildly as only for $\Bbb P_1$.  See more
 in \cite[\S2,\S3]{Sh:915}.
\endremark
\bigskip

\head {\S2 Remarks on $\pi$-bases} \endhead  \resetall \sectno=2
 \spuriousreset
\bigskip

\definition{\stag{p.2} Definition}  1) ${\Cal A}$ is a $\pi$-base
\ub{if}:
\mr
\item "{$(a)$}"  ${\Cal A} \subseteq [\omega]^{\aleph_0}$
\sn
\item "{$(b)$}"  for some ultrafilter $D$ on $\omega,{\Cal A}$ is a
$\pi$-base of $D$, see below, note that $D$ is necessarily non-principal
\ermn
1A) We say ${\Cal A}$ is a $\pi$-base of $D$ if 
$(\forall B \in D)(\exists A \in {\Cal A})(A \subseteq^* B)$.
\nl
1B) $\pi \chi(D) = \text{ Min}\{|{\Cal A}|:{\Cal A}$ is a 
$\pi$-base of $D\}$. 
\nl 
2) ${\Cal A}$ is a strict $\pi$-base \ub{if}:
\mr
\item "{$(a)$}"  ${\Cal A}$ is a $\pi$-base of some $D$
\sn
\item "{$(b)$}"  no subset of ${\Cal A}$ of cardinality $< |{\Cal A}|$
is a $\pi$-base.
\ermn
3) $D$ has a strict $\pi$-base \ub{when} $D$ has a $\pi$-base 
${\Cal A}$ which is a strict $\pi$-base. 
\nl
4) Sp$^*_{\pi\chi} = \{|{\Cal A}|$: there is a non-principal
ultrafilter $D$ on $\omega$ such that ${\Cal A}$ is a strict
$\pi$-base of $D\}$.
\enddefinition
\bigskip

\definition{\stag{p.4} Definition}  For ${\Cal A} \subseteq
[\omega]^{\aleph_0}$ let Id$_{\Cal A} = \{B \subseteq \omega$: for some
$n < \omega$ and
partition $\langle B_\ell:\ell < n \rangle$ of $B$ for no $A \in {\Cal
A}$ and $\ell <n$ do we have $A \subseteq^* B_\ell\}$.
\enddefinition
\bigskip

\demo{\stag{p.5} Observation}  For ${\Cal A} \subseteq
[\omega]^{\aleph_0}$ we have:
\mr
\item "{$(a)$}"  Id$_{\Cal A}$ is an ideal on ${\Cal P}(\omega)$
including the finite sets, though may be equal to ${\Cal P}(\omega)$
\sn
\item "{$(b)$}" if $B \subseteq \omega$ then:
 $B \in [\omega]^{\aleph_0} \backslash
\text{ Id}_{\Cal A}$ iff there is a (non-principal) 
ultrafilter $D$ on $\omega$ to
which $B$ belongs and ${\Cal A}$ is a $\pi$-base of $D$
\sn
\item "{$(c)$}"  ${\Cal A}$ is a $\pi$-base iff 
$\omega \notin \text{ Id}_{\Cal A}$.
\endroster
\enddemo
\bigskip

\demo{Proof} 
\mn
\ub{Clause (a)}:  Obvious.
\mn
\ub{Clause (b)}:  
\sn
\ub{The ``if" direction}:  Let $D$ be a non-principal ultrafilter on
$\omega$ such that $B \in D$ and ${\Cal A}$ is a $\pi$-base of
$D$.  Now for any $n < \omega$ and partition $\langle B_\ell:\ell <
n\rangle$ of $B$ as $B \in D$ and $D$ is an ultrafilter clearly there
is $\ell <n$ such that $B_\ell \in D$ hence by Definition
\scite{p.2}(1A) there is $A \in {\Cal A}$ such that $A \subseteq^*
B_\ell$.  By the definition of Id$_{\Cal A}$ it follows that $B \notin
\text{ Id}_{\Cal A}$ but $[\omega]^{< \aleph_0} \subseteq 
\text{ Id}_{\Cal A}$ so we are done.
\sn
\ub{The ``only if" direction}:  So we are assuming $B \notin 
\text{ Id}_{\Cal A}$ so as Id$_{\Cal A}$ is an ideal of ${\Cal
P}(\omega)$ there is an ultrafilter $D$ on $\omega$ disjoint to
Id$_{\Cal A}$ such that $B \in D$.  So if $B' \in D$ then $B'
\subseteq \omega \wedge B' \notin \text{ Id}_{\Cal A}$ hence by the
definition of Id$_{\Cal A}$ it follows that $(\exists A \in {\Cal
A})(A \subseteq^* B')$.  By Definition \scite{p.2}(1A) this
means that ${\Cal A}$ is a $\pi$-base of $D$.
\mn
\ub{Clause (c)}:  Follows from clause (b).  \hfill$\square_{\scite{p.3}}$
\enddemo
\bigskip

\demo{\stag{p.3} Observation}  1) If $D$ is an ultrafilter on $\omega$
\ub{then} $D$ has a $\pi$-base of cardinality $\pi\chi(D)$. 
\nl
2) ${\Cal A}$ is a $\pi$-base \ub{iff} for every $n \in
[1,\omega)$ and partition $\langle B_\ell:\ell < n \rangle$ of
$\omega$ to finitely many sets, for some 
$A \in {\Cal A}$ and $\ell < n$ we have $A
\subseteq^* B_\ell$. 
\nl
3) Min$\{\pi \chi(D):D$ a non-principal ultrafilter on $\omega\} =
\text{ Min}\{|{\Cal A}|:{\Cal A}$ is a $\pi$-base$\} = \text{
Min}\{|{\Cal A}|:{\Cal A}$ is a strict $\pi$-base$\}$.  
\enddemo
\bigskip

\demo{Proof}  1) By the definition.
\nl
2) For the ``only if" direction, assume ${\Cal A}$ is a $\pi$-base of 
$D$ then Id$_{\Cal A} \subseteq {\Cal P}(\omega)
\backslash D$ (see the proof of \scite{p.4}) so $\omega \notin 
\text{ Id}_{\Cal A}$ and we are done.

For the ``if" direction, use \scite{p.4}.  
\nl
3) Easy.  \hfill$\square_{\scite{p.3}}$
\enddemo
\bigskip

\proclaim{\stag{642.2A} Theorem}  In $\bold V^{\Bbb P}$ as in
\scite{642.2}, we have $\{\mu,\lambda\} \subseteq
{ \text{\rm Sp\/}}^*_{\pi \chi}$ and $\kappa_2 \notin 
{ \text{\rm Sp\/}}^*_{\pi \chi}$.
\endproclaim
\bigskip

\demo{Proof}  Similar to the proof of \scite{642.2} but with some
additions.  The main change is in the proof of $\Vdash_{\Bbb
P} ``\lambda \in \text{ Sp}_\chi"$.  The main addition is that
choosing ${\frak t}_\alpha$ by induciton on $\alpha$ we also define
${\Cal A}_\alpha$ such that
\mr
\item "{$\otimes'_1$}"  $(a),(b) \quad$ as in $\otimes_1$
\sn
\item "{${{}}$}"  $(c) \quad$ as in $\otimes_2(c)$ but only if $\alpha
\ne 2$ mod $\omega$ (and $\alpha = \beta +1$)
\sn
\item "{${{}}$}"  $(d) \quad {\underset\tilde {}\to A_\alpha}$ is a
$\Bbb P^{{\frak t}_\alpha}_0$-name of an infinite subset of $\omega$
\sn
\item "{${{}}$}"  $(e) \quad$ if $\alpha \ne 2$ mod $\omega$ then
$\Vdash_{\Bbb P^{{\frak t}_\alpha}} {\underset\tilde {}\to A_\alpha} =
\omega$ (or do not define ${\underset\tilde {}\to A_\alpha}$)
\sn
\item "{${{}}$}"  $(f) \quad$ if $\alpha < \beta$ are $= 2$ mod
$\omega$ then $\Vdash_{\Bbb P^{{\frak t}_\beta}_\mu}
``{\underset\tilde {}\to A_\beta} \subseteq^* {\underset\tilde {}\to
A_\alpha}$
\sn
\item "{${{}}$}"  $(g) \quad$ if $\beta = \alpha +1$ and $\beta = 2$ mod
$\omega$ and $\underset\tilde {}\to B$ is a $\Bbb P^{{\frak
t}_\alpha}_\mu$-name of an
\nl

\hskip25pt infinite subset of $\omega$ then
$\Vdash_{\Bbb P^{{\frak t}_\beta}_\mu} ``\underset\tilde {}\to B
\nsubseteq^* A_\alpha$.
\ermn
This addition requires that we also prove
\mr
\item "{$\circledast_3$}"  if ${\frak s} \in {\frak K}$ and
$\underset\tilde {}\to D$ is a $\Bbb P^{\frak s}_0$-name of a filter on
$\omega$ including all co-finite subsets of $\omega$ (such that
$\emptyset \notin D$) \ub{then} for some $({\frak t},\underset\tilde
{}\to A)$ we have
{\roster
\itemitem{ $(a)$ }  ${\frak s} \le_{\frak K} {\frak t}$
\sn
\itemitem{ $(b)$ }  $\Vdash_{\Bbb P^{\frak t}} ``\underset\tilde {}\to
A$ is an infinite subset of $\omega$
\sn
\itemitem{ $(c)$ }  if $\underset\tilde {}\to B$ is a 
$\Bbb P^{\frak s}$-name of an infinite subset of $\omega$ then
$\Vdash_{\Bbb P^{\frak t}} ``\underset\tilde {}\to B \nsubseteq^*
\underset\tilde {}\to A"$.
\endroster}
\ermn
[Why $\circledast_3$ holds?  Without loss of generality 
$\Vdash_{\Bbb P^{\frak s}_0} 
``\underset\tilde {}\to D$ is an ultrafilter on $\omega$".  We can
find a pair $(\Bbb P',{\underset\tilde {}\to A'})$
\mr
\item "{$(\alpha)$}"  $\Bbb P'$ is a c.c.c. forcing notion
\sn
\item "{$(\beta)$}"  $\Bbb P^{\frak s}_0 \lessdot \Bbb P'$ moreover
$\Bbb P' = \Bbb P^{\frak s}_0 * \Bbb Q(\underset\tilde {}\to D)$
\sn
\item "{$(\gamma)$}"  $|\Bbb P'| \le \lambda$
\sn
\item "{$(\delta)$}"  $\Vdash_{\Bbb P'} ``\underset\tilde {}\to A$ is
an almost intersection of $\underset\tilde {}\to D$
(i.e. $\underset\tilde {}\to A \in [\omega]^{\aleph_0}$ and $(\forall
B \in \underset\tilde {}\to D)(A \subseteq^* B))$
\sn
\item "{$(\varepsilon)$}"  ${\underset\tilde {}\to \eta'} \in
{}^\omega \omega$ is the generic of $\Bbb Q[\underset\tilde {}\to D]$
and ${\underset\tilde {}\to A'} = \text{ Rang}(\eta)$ so both are
$\Bbb P'$-names.
\ermn
Now we define ${\frak t}'$: for $i \le \mu$ we choose $\Bbb P^{{\frak
t}'}_i = \Bbb P^{\frak s}_i *_{\Bbb P^s_0} \Bbb P'$ and we choose
${\underset\tilde {}\to \tau^{{\frak t}'}_i}$ naturally.  Let $\langle 
{\underset\tilde {}\to n_\rho}:\rho \in {}^{\omega >}2\rangle$ be a
$\Bbb P^{{\frak t}'}_0$-name listing the members of ${\utt A}$.

Now we choose ${\frak t}$ such that ${\frak t}' \le_{\frak K} {\frak
t}$ and for some $\Bbb P^{\frak t}_0$-name $\underset\tilde {}\to
\rho$ of a member of ${}^\omega 2$ we have $\Vdash_{\Bbb P_{\frak t}} 
``\underset\tilde {}\to \rho \ne \underset\tilde {}\to \nu"$ for any
$\Bbb P_{{\frak t}'}$-name (clearly exists, e.g. when $({\frak
t},{\frak t}')$ is like $({\frak t}',{\frak s})$ above).  Now
$\underset\tilde {}\to A := \{{\utt n_{{\utt \rho} \restriction k}}:k
< \omega\}$ is forced to be an infinite subset of ${\underset\tilde {}\to A'}$,
and if it includes a member of ${\Cal P}(\omega)^{\bold V[\Bbb
P_{\frak s}]}$ or even ${\Cal P}(\omega)^{\bold V[\Bbb P_{\frak t}]}$
we get that $\underset\tilde {}\to \rho$ is from $({}^\omega 2)^{\bold
V[\Bbb P_{{\frak t}'}]}$, contradiction.]
\mr
\item "{$(*)_1$}"  $\mu \in \text{ Sp}^*_{\pi\chi}$, in $\bold V^{\Bbb
P}$, of course.
\ermn
[Why?   As there is a $\subseteq^*$-decreasing sequence $\langle
B_\alpha:\alpha < \mu\rangle$ of sets which generates a (non-principle
ultrafilter).  We can use $B_\alpha$ is the generic of $\Bbb P^{{\frak
t}_\lambda}_{\alpha +1}/P^{{\frak t}_\lambda}_\alpha$.]
\mr 
\item "{$(*)_2$}"  $\kappa_2 \notin \text{ Sp}^*_{\pi \chi}$.
\ermn
[Why?  Toward contradiction assume $p^* \in \Bbb P$ and $p^*
\Vdash_{\Bbb P} ``\underset\tilde {}\to D$ is a non-principal
ultrafilter on $\omega$ and $\{{\underset\tilde {}\to
{\Cal U}_\varepsilon}:\varepsilon < \kappa_2\}$ is a sequence of infinite
subsets of $\omega$ which is a strict $\pi$-base of $\underset\tilde
{}\to D"$; so $p^* \Vdash_{\Bbb P} ``\{{\underset\tilde {}\to
{\Cal U}_\varepsilon}:\varepsilon < \zeta\}$ is not a $\pi$-base of any
ultrafilter on $\omega$" for every $\zeta < \kappa_2$, hence for some
$\langle {\underset\tilde {}\to B_{\zeta,\ell}}:\ell < 
{\underset\tilde {}\to n_\zeta}\rangle$ we have $p^* \Vdash 
``{\underset\tilde {}\to n_\ell} < \omega$ and $\langle
{\underset\tilde {}\to B_{\zeta,\ell}}:\ell < {\underset\tilde {}\to
n_\ell}\rangle$ is a partition of $\omega$ and $\varepsilon < \zeta
\wedge \ell < {\underset\tilde {}\to n_\zeta} \Rightarrow 
{\Cal U}_\varepsilon \nsubseteq^* {\underset\tilde {}\to B_{\zeta,\ell}}"$.
We now as in the proof of \scite{642.2}, choose suitable $\beta <
\lambda$ and consider $\langle {\underset\tilde {}\to B^*_\ell}:\ell <
\underset\tilde {}\to n\rangle = \hat{\bold
j}^{-1}_\beta(\langle{\underset\tilde {}\to B_{\zeta,\ell}}:\ell <
{\underset\tilde {}\to n_\zeta}\rangle:\zeta < \kappa_2\rangle/{\Cal
D}_2)$ so $p^* \Vdash_{\Bbb P^{\beta +1}} ``\langle {\underset\tilde
{}\to B^*_\ell}:\ell < \underset\tilde {}\to n\rangle$ is a partition
of $\omega$ to finitely many sets and $\varepsilon < \kappa_2 \wedge
\ell < \underset\tilde {}\to n \Rightarrow {\underset\tilde {}\to
{\Cal U}_\varepsilon} \nsubseteq^* {\underset\tilde {}\to B^*_\ell}"$.  But
this contradicts $p^* \Vdash_{\Bbb P} ``\{{\underset\tilde {}\to
{\Cal U}_\varepsilon}:\varepsilon < \kappa_2\}$ is a $\pi$-base.]
\mr
\item "{$(*)_3$}"  $\lambda \in \text{ Sp}^*_{\pi}$.
\ermn
[Why?  Clearly it is forced (i.e. $\Vdash_{\Bbb P_\lambda}$) that
$\langle {\underset\tilde {}\to A_{\omega \alpha +2}}:\alpha <
\lambda\rangle$ is a $\subseteq^*$-decreasing sequence of infinite
subsets of $\omega$, hence there is an ultrafilter of $D$ on $\omega$
including it.  Now ${\underset\tilde {}\to A_{\omega \alpha +2}}$
witness that 
${\Cal P}(\omega)^{\bold V[\Bbb P_{{\frak t}_{\omega \alpha +2}}]}$ is
not a $\pi$-base of $\underset\tilde {}\to D$ (recalling clause (h) of
$\circledast'_1$).  As $\lambda$ is regular we are done.]
\hfill$\square_{\scite{642.2A}}$
\enddemo

\bigskip\bigskip
    
REFERENCES.  
\bibliographystyle{lit-plain}
\bibliography{lista,listb,listx,listf,liste}

\enddocument